\documentclass{article}
\usepackage{amssymb}
\usepackage{amsmath}
\usepackage[all]{xy}

\def\Im{\mathop{\rm Im}\nolimits}
\def\Ker{\mathop{\rm Ker}\nolimits}
\def\Coker{\mathop{\rm Coker}\nolimits}
\def\lfd{\mathop{\rm l.fd}\nolimits}
\def\lpd{\mathop{\rm l.pd}\nolimits}
\def\lid{\mathop{\rm l.id}\nolimits}
\def\rfd{\mathop{\rm r.fd}\nolimits}
\def\rpd{\mathop{\rm r.pd}\nolimits}
\def\rid{\mathop{\rm r.id}\nolimits}
\def\lfindim{\mathop{\rm l.fin.dim}\nolimits}
\def\rfindim{\mathop{\rm r.fin.dim}\nolimits}
\def\domdim{\mathop{\rm dom.dim}\nolimits}
\def\grade{\mathop{\rm grade}\nolimits}
\def\sgrade{\mathop{\rm s.grade}\nolimits}
\def\mod{\mathop{\rm mod}\nolimits}

\title{\Large \bf On the Grade of Modules over Noetherian Rings
}
\author{Zhaoyong Huang\thanks{{\it E-mail address}: huangzy@nju.edu.cn}
\\{\small \it Department of Mathematics, Nanjing University,
Nanjing 210093, P. R. China}}
\usepackage{amssymb}
\date{}
\begin{document}
\baselineskip=18pt \maketitle

\begin{abstract} Let $\Lambda$ be a left and right noetherian ring
and $\mod \Lambda$ the category of finitely generated left
$\Lambda$-modules. In this paper we show the following results: (1)
For a positive integer $k$, the condition that the subcategory of
$\mod \Lambda$ consisting of $i$-torsionfree modules coincides with
the subcategory of $\mod \Lambda$ consisting of $i$-syzygy modules
for any $1\leq i \leq k$ is left-right symmetric. (2) If $\Lambda$
is an Auslander ring and $N$ is in $\mod \Lambda ^{op}$ with $\grade
N=k<\infty$, then $N$ is pure of grade $k$ if and only if $N$ can be
embedded into a finite direct sum of copies of the $(k+1)$st term in
a minimal injective resolution of $\Lambda$ as a right
$\Lambda$-module. (3) Assume that both the left and right
self-injective dimensions of $\Lambda$ are $k$. If $\grade {\rm
Ext}_{\Lambda}^k(M, \Lambda)\geq k$ for any $M\in\mod \Lambda$ and
$\grade {\rm Ext}_{\Lambda}^i(N, \Lambda)\geq i$ for any $N\in\mod
\Lambda ^{op}$ and $1\leq i \leq k-1$, then the socle of the last
term in a minimal injective resolution of $\Lambda$ as a right
$\Lambda$-module is non-zero.
\end{abstract}

\vspace{0.5cm}

\centerline{\large \bf 1. Introduction}

\vspace{0.2cm}

Throughout this paper, $\Lambda$ is a left and right noetherian ring
and $\mod \Lambda$ denotes the category of finitely generated left
$\Lambda$-modules. It is well known that the properties of grade of
modules are useful to characterize rings as well as to study the
dual properties of modules (see, for example, [2], [4], [7], [10]
and [13--15]). In this paper, we study the homological properties of
modules over noetherian rings, especially over $k$-Gorenstein rings
and related rings, under some grade conditions of modules.

Let $M$ be in $\mod \Lambda$ (resp. $\mod \Lambda ^{op}$). For a
non-negative integer $t$, we say that the {\it grade} of $M$ is
equal to $t$, denoted by $\grade M=t$, if Ext$_{\Lambda}^i(M,
\Lambda)=0$ for any $0\leq i<t$ and Ext$_{\Lambda}^t(M, \Lambda)\neq
0$. We say the {\it strong grade} of $M$ is equal to $t$, denoted by
$\sgrade M=t$, if $\grade A=t$ for each non-zero submodule $A$ of
$M$. Moreover, if Ext$_{\Lambda}^i(M, \Lambda)=0$ for any $i \geq
0$, then we write $\grade M=\infty$ (see [2]).

For a positive integer $k$, Auslander and Bridger introduced in [2]
the notion $k$-torsionfree modules. Such a class of modules is
natural and interesting in homological algebra. It was showed in [2,
Theorem 2.17] that a $k$-torsionfree module is $k$-syzygy; but the
converse is not true in general. It was showed in [2, Proposition
2.26] that each $i$-syzygy module is $i$-torsionfree for any $1 \leq
i \leq k$ if and only if $\grade {\rm Ext}_{\Lambda}^{i+1}(M,
\Lambda)\geq i$ for any $M\in\mod \Lambda$ and $1 \leq i \leq k-1$.
We show in Section 2 that this result is left-right symmetric. Under
these equivalent conditions we show that $\lfindim \Lambda \leq k$
if and only if $\Lambda$ satisfies the condition that $N=0$ for
every $N \in \Lambda ^{op}$ satisfying $\grade N \geq k+1$. As a
corollary, we get that if the left flat dimension of the $i$th term
in a minimal injective resolution of $\Lambda$ as a left
$\Lambda$-module is at most $i$ for all $i$, then the left
self-injective dimension of $\Lambda$ and its small left finitistic
dimension are identical, and the difference between the right
self-injective dimension of $\Lambda$ and its small right finitistic
dimension is at most one.

Non-commutative Gorenstein rings are already defined. These are
rings for which the (left/right) self-injective dimension of the
ring is finite. In addition, recall that $\Lambda$ is called
$k$-{\it Gorenstein} if the right flat dimension of the $i$th term
in a minimal injective resolution of $\Lambda$ as a right
$\Lambda$-module is at most $i-1$ for any $1\leq i \leq k$; and
$\Lambda$ is called an {\it Auslander ring} if it is $k$-Gorenstein
for all $k$ (see [10] and [11]). Auslander gave some useful
equivalent conditions of $k$-Gorenstein rings in term of the right
flat dimension and grade of modules as follows, which shows that the
notion of $k$-Gorenstein rings is left-right symmetric.

\vspace{0.2cm}

{\bf Auslander's Theorem} ([10, Theorem 3.7]) {\it The following
statements are equivalent for} $\Lambda$.

(1) $\Lambda$ {\it is a} $k$-{\it Gorenstein ring}.

(2) {\it The left flat dimension of the} $i$th {\it term in a
minimal injective resolution of} $\Lambda$ {\it as a left}
$\Lambda$-{\it module is at most} $i-1$ {\it for any} $1\leq i \leq
k$.

(3) $\sgrade {\rm Ext}_{\Lambda}^i(M, \Lambda)\geq i$ {\it for any}
$M \in\mod \Lambda$ {\it and} $1\leq i \leq k$.

(4) $\sgrade {\rm Ext}_{\Lambda}^i(N, \Lambda)\geq i$ {\it for any}
$N \in\mod \Lambda ^{op}$ {\it and} $1\leq i \leq k$.

\vspace{0.2cm}

In Section 3 we study the purity of modules over $k$-Gorenstein
rings. Let $k$ be a non-negative integer and $\Lambda$ a
$(k+1)$-Gorenstein ring. For a module $M\in\mod \Lambda ^{op}$ with
$\grade N=k<\infty$, we show that $N$ is pure of grade $k$ if and
only if $N$ can be embedded into a finite direct sum of copies of
the $(k+1)$st term in a minimal injective resolution of $\Lambda$ as
a right $\Lambda$-module.

Assume that both the left and right self-injective dimensions of
$\Lambda$ are $k$. If $\grade {\rm Ext}_{\Lambda}^k(M, \Lambda)\geq
k$ for any $M\in\mod \Lambda$ and $\grade {\rm Ext}_{\Lambda}^i(N,
\Lambda)\geq i$ for any $N\in\mod \Lambda ^{op}$ and $1\leq i \leq
k-1$, we show in Section 4 that the socle of the last term in a
minimal injective resolution of $\Lambda$ as a right
$\Lambda$-module is non-zero. As an immediate result, we have that
if $\Lambda$ is $(k-1)$-Gorenstein with both the left and right
self-injective dimensions $k$, then the socle of the last term in a
minimal injective resolution of $\Lambda$ as a right
$\Lambda$-module is also non-zero. Some known results are obtained
as corollaries.

In the following, we assume that
$$0\to \Lambda \to I_0 \to I_1 \to \cdots \to I_i \to \cdots$$
is a minimal injective resolution of $\Lambda$ as a right
$\Lambda$-module, and that
$$0\to \Lambda \to I'_0 \to I'_1 \to \cdots \to I'_i \to \cdots$$
is a minimal injective resolution of $\Lambda$ as a left
$\Lambda$-module. For a left $\Lambda$-module $M$,
$\lfd_{\Lambda}M$, $\lpd_{\Lambda}M$ and $\lid_{\Lambda}M$ denote
the left flat dimension, the left projective dimension and the left
injective dimension of $M$ respectively, and for a right
$\Lambda$-module $N$, $\rfd_{\Lambda}N$, $\rpd_{\Lambda}N$ and
$\rid_{\Lambda}N$ denote the right flat dimension, the right
projective dimension and the right injective dimension of $N$
respectively.

\vspace{0.5cm}

\centerline{\large \bf 2. $k$-torsionfree modules and $k$-syzygy
modules}

\vspace{0.2cm}

Let $M$ be in $\mod \Lambda$ (resp. $\mod \Lambda ^{op}$), and let
$\sigma _{M}: M \to M^{**}$ via $\sigma _{M}(x)(f)=f(x)$ for any
$x\in M$ and $f\in M^*$ be the canonical evaluation homomorphism,
where $(\ )^*=$Hom$_{\Lambda}(-, \Lambda)$. $M$ is called a {\it
torsionless module} if $\sigma _{M}$ is a monomorphism; and $M$ is
called a {\it reflexive module} if $\sigma _{M}$ is an isomorphism.
For a positive integer $k$, we call $M$ a $k$-{\it syzygy} module if
there is an exact sequence
$$0 \to M \to Q_{0}
\to Q_{1} \to \cdots \to Q_{k-1}$$ with all $Q_{i}$ finitely
generated projective. On the other hand, assume that
$$P_{1} \buildrel {f} \over \longrightarrow P_{0} \to M \to 0$$ is a
resolution of $M$ with $P_0$ and $P_1$ finitely generated
projective. Then we get an exact sequence
$$0 \to M^* \to P_{0}^*
\buildrel {f^*} \over \longrightarrow P_{1}^* \to X \to 0,$$ where
$X=\Coker f^*$. $M$ is called a $k$-{\it torsionfree module} if
Ext$_{\Lambda}^i(X, \Lambda)=0$ for any $1\leq i \leq k$ (see [2]).
Because $X$ is unique up to projective summands, the notion of
$k$-torsionfree modules is well defined. By [12, Lemma 1.5], we have
the exact sequence
$$0\to {\rm Ext}_{\Lambda}^1(X, \Lambda) \to M \buildrel {\sigma
_M} \over \longrightarrow M^{**} \to {\rm Ext}_{\Lambda}^{2}(X,
\Lambda) \to 0.$$ So $M$ is 1-torsionfree if and only if it is
torsionless if and only if it is 1-syzygy; and that $M$ is
2-torsionfree if and only if it is reflexive. We use
$\mathcal{T}^k(\mod \Lambda)$ (resp. $\Omega ^k(\mod \Lambda)$) to
denote the full subcategory of $\mod \Lambda$ consisting of
$k$-torsionfree modules (resp. $k$-syzygy modules). It is known from
[2, Theorem 2.17] that $\mathcal{T}^k(\mod \Lambda) \subseteq \Omega
^k(\mod \Lambda)$. So, it is natural to ask when they are identical.
Auslander and Bridger in [2] gave the following result.

\vspace{0.2cm}

{\bf Proposition 2.1} ([2, Proposition 2.26] or [4, Proposition
1.6]) {\it For a positive integer} $k$, {\it the following
statements are equivalent}.

(1) $\grade {\rm Ext}_{\Lambda}^{i+1}(M, \Lambda)\geq i$ {\it for
any} $M\in\mod \Lambda$ {\it and} $1\leq i \leq k-1$.

(2) $\Omega ^i(\mod \Lambda)=\mathcal{T}^i(\mod \Lambda)$ {\it for
any} $1\leq i \leq k$.

\vspace{0.2cm}

In this section we show this result is left-right symmetric. To
get our result we need two lemmas.

\vspace{0.2cm}

{\bf Lemma 2.2} ([12, Lemma 1.6]) {\it The following statements are
equivalent}.

(1) $[{\rm Ext}_{\Lambda}^2(M, \Lambda)]^*=0$ {\it for any}
$M\in\mod \Lambda$.

(2) $[{\rm Ext}_{\Lambda}^2(N, \Lambda)]^*=0$ {\it for any}
$N\in\mod \Lambda ^{op}$.

(3) $M^*$ {\it is reflexive for any} $M\in\mod \Lambda$.

(4) $N^*$ {\it is reflexive for any} $N\in\mod \Lambda ^{op}$.

\vspace{0.2cm}

{\bf Lemma 2.3} ([12, Lemma 6.2]) {\it Let} $n$ {\it be a
non-negative integer and} $X\in\mod \Lambda$ ({\it resp}. $\mod
\Lambda ^{op}$). {\it If} $\grade X\geq n$ {\it and} $\grade {\rm
Ext}_{\Lambda}^n(X, \Lambda)\geq n+1$, {\it then}
Ext$_{\Lambda}^n(X, \Lambda)=0$.

\vspace{0.2cm}

{\bf Theorem 2.4} {\it For a positive integer} $k$, {\it the
following statements are equivalent}.

(1) $\grade {\rm Ext}_{\Lambda}^{i+1}(M, \Lambda)\geq i$ {\it for
any} $M\in\mod \Lambda$ {\it and} $1\leq i \leq k-1$.

(2) $\Omega ^i(\mod \Lambda)=\mathcal{T}^i(\mod \Lambda)$ {\it for
any} $1\leq i \leq k$.

(3) $\grade {\rm Ext}_{\Lambda}^{i+1}(N, \Lambda)\geq i$ {\it for
any} $N\in\mod \Lambda ^{op}$ {\it and} $1\leq i \leq k-1$.

(4) $\Omega ^i(\mod \Lambda ^{op})=\mathcal{T}^i(\mod \Lambda
^{op})$ {\it for any} $1\leq i \leq k$.

\vspace{0.2cm}

{\it Proof.} By Proposition 2.1 and its dual statement, we get the
equivalence of (1) and (2) and that of (3) and (4). In the following
we prove (3) implies (1) by induction on $k$. The case $k=1$ is
trivial. The case $k=2$ follows from Lemma 2.2. Now suppose $k\geq
3$.

Let $M\in\mod \Lambda$ and
$$\cdots \to P_i \to \cdots \to P_1 \to P_0 \to M \to 0$$ a
projective resolution of $M$ in $\mod \Lambda$. Put $M_i=\Coker (P_i
\to P_{i-1})$ (where $M_1=M$) and $X_i=\Coker (P_{i-1}^* \to P_i^*)$
for any $i\geq 1$. By the induction assumption, we have $\grade {\rm
Ext}_{\Lambda}^{i+1}(M, \Lambda)\geq i$ for any $1\leq i \leq k-2$
and $\grade {\rm Ext}_{\Lambda}^k(M, \Lambda)\geq k-2$. So it
suffices to prove Ext$_{\Lambda}^{k-2}$(Ext$_{\Lambda}^k(M,
\Lambda), \Lambda)=0$.

By Proposition 2.1, $\Omega ^i(\mod \Lambda)=\mathcal{T}^i(\mod
\Lambda)$ for any $1\leq i \leq k-1$. Since $M_t \in \Omega
^{k-1}(\mod \Lambda)$ for any $t\geq k$, $M_t \in
\mathcal{T}^{k-1}(\mod \Lambda)$ for any $t\geq k$. It follows that
Ext$_{\Lambda}^i(X_t, \Lambda)=0$ for any $1 \leq i \leq k-1$ and
$t\geq k$.

On the other hand, by [14, Lemma 2] we have the exact sequence
$$0\to {\rm Ext}_{\Lambda}^k(M, \Lambda) \to X_k \to P_{k+1}^* \to
X_{k+1} \to 0.$$ Put $K=\Im (X_k \to P_{k+1}^*)$. From the exactness
of $0\to K \to P_{k+1}^* \to X_{k+1} \to 0$ we know that
Ext$_{\Lambda}^i(K, \Lambda)=0$ for any $1 \leq i \leq k-2$ and
${\rm Ext}_{\Lambda}^{k-1}(K, \Lambda)\cong {\rm
Ext}_{\Lambda}^k(X_{k+1}, \Lambda)$. Moreover, from the exactness of
$0\to {\rm Ext}_{\Lambda}^k(M, \Lambda) \to X_k \to K \to 0$ we know
that ${\rm Ext}_{\Lambda}^{k-1}(K, \Lambda)\cong {\rm
Ext}_{\Lambda}^{k-2}({\rm Ext}_{\Lambda}^k(M, \Lambda), \Lambda)$.
So ${\rm Ext}_{\Lambda}^{k-2}({\rm Ext}_{\Lambda}^k(M, \Lambda),
\Lambda)\cong {\rm Ext}_{\Lambda}^k(X_{k+1}, \Lambda)$. By (3) we
then have $\grade {\rm Ext}_{\Lambda}^{k-2}({\rm Ext}_{\Lambda}^k(M,
\Lambda), \Lambda)=\grade {\rm Ext}_{\Lambda}^k(X_{k+1},
\Lambda)\geq k-1$. It follows from Lemma 2.3 that
Ext$_{\Lambda}^{k-2}$(Ext$_{\Lambda}^k(M, \Lambda), \Lambda)=0$.
Dually we have (1) implies (3). $\blacksquare$

\vspace{0.2cm}

{\bf Corollary 2.5} {\it The following statements are equivalent}.

(1) $\grade {\rm Ext}_{\Lambda}^{i+1}(M, \Lambda)\geq i$ {\it for
any} $M\in\mod \Lambda$ {\it and} $i \geq 1$.

(2) $\Omega ^i(\mod \Lambda)=\mathcal{T}^i(\mod \Lambda)$ {\it for
any} $i \geq 1$.

(3) $\grade {\rm Ext}_{\Lambda}^{i+1}(N, \Lambda)\geq i$ {\it for
any} $N\in\mod \Lambda ^{op}$ {\it and} $i \geq 1$.

(4) $\Omega ^i(\mod \Lambda ^{op})=\mathcal{T}^i(\mod \Lambda
^{op})$ {\it for any} $i \geq 1$.

\vspace{0.2cm}

{\it Remark.} (1) We remark that if $\rfd_{\Lambda}I_i \leq i+1$ for
any $0\leq i \leq k-2$ (especially, if $\Lambda$ is
$(k-1)$-Gorenstein) then the condition (1) in Theorem 2.4 is
satisfied by [4, Proposition 2.2].

(2) Consider the following grade conditions.

$(a_k)$ $\sgrade {\rm Ext}_{\Lambda}^i(M, \Lambda)\geq i$ for any $M
\in\mod \Lambda$ and $1\leq i \leq k$.

$(b_k)$ $\grade {\rm Ext}_{\Lambda}^i(M, \Lambda)\geq i$ for any $M
\in\mod \Lambda$ and $1\leq i \leq k$.

$(c_k)$ $\sgrade {\rm Ext}_{\Lambda}^{i+1}(M, \Lambda)\geq i$ for
any $M \in\mod \Lambda$ and $1\leq i \leq k$.

$(d_k)$ $\grade {\rm Ext}_{\Lambda}^{i+1}(M, \Lambda)\geq i$ for any
$M \in\mod \Lambda$ and $1\leq i \leq k$.

About the symmetry of these conditions, the following facts are
known: Both $(a_k)$ and $(d_k)$ are left-right symmetric by
Auslander's Theorem and Theorem 2.4, respectively. Neither $(b_k)$
nor $(c_k)$ is left-right symmetric by [13, Theorem 3.3 and p.1460
``Remark"]. On the other hand, we clearly have the implications
$(a_k)\Rightarrow (b_k)\Rightarrow (d_k)$ and $(a_k)\Rightarrow
(c_k)\Rightarrow (d_k)$. However, by the above argument, none of the
implications cannot be reversed.

\vspace{0.2cm}

{\bf Corollary 2.6} {\it Let} $\Lambda$ {\it be an artinian ring. If
one of the equivalent conditions in Theorem 2.4 is satisfied for}
$\Lambda$, {\it then} $\lid_{\Lambda}\Lambda \leq k$ {\it if and
only if} $\rid_{\Lambda}\Lambda \leq k$. {\it In particular, if one
of the equivalent conditions in Corollary 2.5 is satisfied for}
$\Lambda$, {\it then} $\lid_{\Lambda}\Lambda=\rid_{\Lambda}\Lambda$.

\vspace{0.2cm}

{\it Proof.} By assumption we have $\grade {\rm Ext}_{\Lambda}^k(M,
\Lambda)\geq k-1$ for any $M\in\mod \Lambda$. If
$\rid_{\Lambda}\Lambda \leq k$, then $\lid_{\Lambda}\Lambda \leq
2k-1$ by [14, Theorem]. It follows from [19, Lemma A] that
$\lid_{\Lambda}\Lambda=\rid_{\Lambda}\Lambda (\leq k)$. The converse
is proved dually. The latter assertion follows easily from the
former one. $\blacksquare$

\vspace{0.2cm}

Recall that the {\it small left finitistic dimension} of $\Lambda$,
written $\lfindim \Lambda$, is defined to be
sup$\{\lpd_{\Lambda}M|M$ is in $\mod \Lambda$ with
$\lpd_{\Lambda}M<\infty \}$. The {\it small right finitistic
dimension} of $\Lambda$ is defined dually, and is denoted by
$\rfindim \Lambda$. As an application of Theorem 2.4, we will prove
the following

\vspace{0.2cm}

{\bf Theorem 2.7} {\it Let} $k$ {\it be a non-negative integer. If
one of the equivalent conditions in Theorem 2.4 is satisfied for}
$\Lambda$, {\it then the
following statements are equivalent.}

(1) $\lfindim \Lambda \leq k$.

(2) $N=0$ {\it for every} $N\in\mod \Lambda ^{op}$ {\it satisfying}
$\grade N \geq k+1$.

\vspace{0.2cm}

To prove this theorem, we need some lemmas.

Let $n$ be a positive integer and
$$P_n \buildrel {d_n} \over \longrightarrow P_{n-1}
\buildrel {d_{n-1}} \over \longrightarrow \cdots \buildrel {d_2}
\over \longrightarrow P_2 \buildrel {d_1} \over \longrightarrow P_1
\to A \to 0$$ a projective resolution of $A$ in $\mod \Lambda$
(resp. $\mod \Lambda ^{op}$). Put $X=\Coker d_n^*$.

\vspace{0.2cm}

{\bf Lemma 2.8} ([12, Lemma 1.5]) {\it Let} $A$ {\it and} $X$ {\it
be as above. If} Ext$_{\Lambda}^i(A, \Lambda)=0$ {\it for any}
$1\leq i \leq n-1$, {\it then we have the following exact sequence}:
$$0\to {\rm Ext}_{\Lambda}^n(X, \Lambda) \to A \buildrel {\sigma
_A} \over \longrightarrow A^{**} \to {\rm Ext}_{\Lambda}^{n+1}(X,
\Lambda) \to 0.$$ {\it In particular, when} $n=1$, {\it we obtain
the exact sequence}
$$0\to {\rm Ext}_{\Lambda}^1(X, \Lambda) \to A \buildrel {\sigma
_A} \over \longrightarrow A^{**} \to {\rm Ext}_{\Lambda}^2(X,
\Lambda) \to 0,$$ {\it where} $X=\Coker d_1^*$.

\vspace{0.2cm}

{\bf Lemma 2.9} $\lfindim \Lambda=0$ {\it if and only if} $N=0$ {\it
for every} $N\in\mod \Lambda ^{op}$ {\it satisfying} $N^*=0$.

\vspace{0.2cm}

{\it Proof.} By [5, Corollary 5.6 and Theorem 5.4]. $\blacksquare$

\vspace{0.2cm}

{\bf Lemma 2.10} $\lfindim \Lambda \leq 1$ {\it if and only if}
$N=0$ {\it for every} $N\in\mod \Lambda ^{op}$ {\it satisfying}
$N^*=0=$Ext$_{\Lambda}^1(N, \Lambda)$.

\vspace{0.2cm}

{\it Proof.} {\it The necessity}. Let $N$ be in $\mod \Lambda ^{op}$
with $N^*=0=$Ext$_{\Lambda}^1(N, \Lambda)$ and
$$0\to K \to P \to N \to 0$$ an exact sequence in $\mod \Lambda ^{op}$
with $P$ projective. Then $K^* (\cong P^*)$ is projective. By [5,
Proposition 5.3], $K$ is projective. It then follows from Lemma 2.8
that we have an exact sequence:
$$0\to {\rm Ext}_{\Lambda}^1({\rm Ext}_{\Lambda}^1(N, \Lambda), \Lambda)
\to N \buildrel {\sigma _N} \over \to N^{**}.$$
But $N^*=0=$Ext$_{\Lambda}^1(N, \Lambda)$, so $N=0$.

{\it The sufficiency}. It suffices to show that if there is an exact
sequence:
$$0 \to P_1 \to P_0 \to M \to 0$$ in $\mod \Lambda$ with $P_0$ and
$P_1$ projective and $M$ torsionless, then $M$ is projective. Put
$N=\Coker (P_0^* \to P_1^*)$. Then $N^*=0$ since $P_0$ and $P_1$ are
reflexive. By Lemma 2.8 we have that Ext$_{\Lambda}^1(N, \Lambda)
\cong \Ker \sigma _M =0$. Then $N=0$ by assumption. So $M^*$ and
$M^{**}$ are projective. On the other hand, we have an exact
sequence
$$0 \to P_1^{**} \to P_0^{**} \to M^{**} \to 0.$$ Thus $M$ is reflexive
and therefore it is projective. $\blacksquare$

\vspace{0.2cm}

The next lemma finishes the proof of Theorem 2.7 in one direction.

\vspace{0.2cm}

{\bf Lemma 2.11} {\it Let} $k$ {\it be a non-negative integer. If}
$\lfindim \Lambda \leq k$, {\it then} $N=0$ {\it for every}
$N\in\mod \Lambda ^{op}$ {\it satisfying} $\grade N \geq k+1$.

\vspace{0.2cm}

{\it Proof.} The cases $k=0$ and $k=1$ were proved in Lemmas 2.9 and
2.10 respectively. Now suppose $k\geq 2$.

Let $N$ be in $\mod \Lambda ^{op}$ with $\grade N \geq k+1$ and
$$\cdots \to Q_{k+1} \to Q_k \to \cdots \to Q_1 \to Q_0 \to N \to 0$$
a projective resolution of $N$ in $\mod \Lambda ^{op}$. Put
$Y=\Coker (Q_k^* \to Q_{k+1}^*)$. Then $\lpd_{\Lambda}Y\leq k+1$. By
assumption $\lfindim \Lambda \leq k$, hence $\lpd_{\Lambda}Y\leq k$.
On the other hand, by Lemma 2.8 we have that
$N\cong$Ext$_{\Lambda}^{k+1}(Y, \Lambda)$. It follows that $N=0$.
$\blacksquare$

\vspace{0.2cm}

We now recall some notions from [17].  A monomorphism $X^{**}
\buildrel {\rho ^*}\over \to Y^*$ in $\mod \Lambda$ (resp. $\mod
\Lambda ^{op}$) is called a {\it double dual embedding} if it is the
dual of an epimorphism $Y \buildrel {\rho}\over \to X^*$ in $\mod
\Lambda ^{op}$ (resp. $\mod \Lambda$). For a positive integer $k$, a
torsionless module $T_k$ in $\mod \Lambda$ (resp. $\mod \Lambda
^{op}$) is said to be of $D$-{\it class} $k$ if it can be fitted
into a diagram of the form

$$\xymatrix{& & & & 0\ar[r] & T_{k-1}^{**}\ar[r]
& P_{k-1}\ar[r] & T_k \ar[r] & 0 \\
& & & \cdots \ar[r] & P_{k-2}\ar[r] & T_{k-1}\ar[r]\ar[u]_{\sigma
_{T_{k-1}}} & 0 & & \\
& & & \cdots & & & & & \\
& & 0\ar[r] & T_{2}^{**}\ar[r] & \cdots & &
& & \\
0 \ar[r] & T_1^{**} \ar[r] & P_1 \ar[r] & T_2 \ar[r]\ar[u]_{\sigma
_{T_{2}}} & 0 & & & & }
$$ where each $P_i$ is projective in $\mod \Lambda$ (resp. $\mod
\Lambda ^{op}$) and the horizontal monomorphisms are double dual
embeddings. Any torsionless module is said to be of $D$-{\it class
1}.

\vspace{0.2cm}

{\bf Lemma 2.12} ([17, Theorem 4.2]) {\it For a positive integer}
$k$, {\it the following statements are equivalent}.

(1) $\lfindim \Lambda \leq k$.

(2) {\it The only modules of} $D$-{\it class} $k$ {\it in} $\mod
\Lambda ^{op}$ {\it with projective duals are the projective
modules}.

\vspace{0.2cm}

{\bf Lemma 2.13} {\it Assume that} $\Lambda$ {\it satisfies one of
the equivalent conditions in Theorem 2.4 and let} $k$ {\it be a
non-negative integer. If} $\Lambda$ {\it in addition satisfies the
condition that} $N=0$ {\it for every} $N\in\mod \Lambda ^{op}$ {\it
satisfying} $\grade N \geq k+1$, {\it then} $\lfindim \Lambda \leq
k$.

\vspace{0.2cm}

{\it Proof.} By Lemmas 2.9 and 2.10, we only need to prove the case
$k\geq 2$.

Let $T_k$ be a module of $D$-class $k$ in $\mod \Lambda ^{op}$ with
$T_k^*$ projective (where $k\geq 2$). From the proof of [17, Theorem
2.1] we know that there are exact sequences
$$0 \to T_i^{**} \buildrel {\rho _i^*} \over \to
P_i^* \buildrel {\pi _i^*} \over \to T_{i+1} \to 0
\eqno{(i)^*}$$
$$0 \to T_{i+1}^* \buildrel {\sigma _{P_i}^{-1}\pi _i^*}
\over \longrightarrow P_i \buildrel {\rho _i} \over \to T_i^* \to 0
\eqno{(i)}$$ for any $1 \leq i \leq k-1$, where all $P_i$ are
projective in $\mod \Lambda$ and all $T_i$ are of $D$-class $i$ in
$\mod \Lambda ^{op}$. Then we have an exact sequence
$$0 \to T_{k-1}^* \to P_{k-2} \to \cdots \to
P_2 \to P_1 \to T_1^* \to 0.$$ It is trivial
that $T_1^*$ is 2-syzygy.
So $T_{k-1}^*$ is $k$-syzygy and hence
it is $k$-torsionfree by Theorem 2.4.

The exact sequence $(k-1)$ induces an exact sequence
$$0 \to T_{k-1}^{**} \buildrel {\rho _{k-1}^*} \over \longrightarrow
P_{k-1}^* \buildrel {\pi _{k-1}^*(\sigma _{P_{k-1}}^{-1})^*} \over
\longrightarrow T_k^{**} \to N \to 0 \eqno{(\dag)}$$ where $N=\Coker
(\pi _{k-1}^*(\sigma _{P_{k-1}}^{-1})^*)$. Since $T_k^*$ is
projective and $T_{k-1}^*$ is $k$-torsionfree, Ext$_{\Lambda}^i(N,
\Lambda)=0$ for any $1\leq i \leq k$. On the other hand, notice that
$T_k^*$ and $P_{k-1}$ are reflexive, it then follows that $N^*=0$
and so $\grade N\geq k+1$. Thus, by assumption, we have that $N=0$
and the exact sequence $(\dag)$ splits.

By [1, Proposition 20.14], we have $(\sigma _{P_{k-1}}^{-1})^*
=(\sigma _{P_{k-1}}^*)^{-1}=\sigma _{P_{k-1}^*}$. We in addition
note that $\sigma _{T_k}\pi _{k-1}=\pi _{k-1}^{**}\sigma
_{P_{k-1}^*}$. So we have the following commutative diagram with
exact rows
$$\xymatrix{0\ar[r] & T_{k-1}^{**}\ar[r]^{\rho _{k-1}^*}\ar@{=}[d]
& P_{k-1}^*\ar[rr]^{\pi _{k-1}}\ar@{=}[d]
& & T_k\ar[r]\ar[d]^{\sigma _{T_k}} & 0\\
0\ar[r] & T_{k-1}^{**}\ar[r]^{\rho _{k-1}^*} & P_{k-1}^*\ar[rr]^{\pi
_{k-1}^{**}(\sigma _{P_{k-1}}^{-1})^*} & & T_k^{**}\ar[r] & 0 }$$
Then it is trivial that $\sigma _{T_k}$ is an isomorphism and so
$T_k$ is projective. It follows from Lemma 2.12 that $\lfindim
\Lambda \leq k$. $\blacksquare$

\vspace{0.2cm}

Now Theorem 2.7 follows from Lemmas 2.11 and 2.13. $\blacksquare$

\vspace{0.2cm}

It was showed in [15, Theorem 2.12] that for a $(k+1)$-Gorenstein
ring $\Lambda$, if $\lfindim \Lambda =k$ then $\lid_{\Lambda}\Lambda
\leq k$. The following corollary generalizes this result.

\vspace{0.2cm}

{\bf Corollary 2.14} {\it Let} $k$ {\it be a non-negative integer.
If} $\lfd_{\Lambda}I'_i\leq i+1$ {\it for any} $0\leq i \leq k$,
{\it then we have}

(1) $\lid_{\Lambda}\Lambda=k$ {\it if and only if} $\lfindim \Lambda
=k$.

(2) $k\leq\rid_{\Lambda}\Lambda \leq k+1$ {\it if} $\rfindim \Lambda
=k$.

{\it In particular, if} $\Lambda$ {\it is} $(k+1)$-{\it Gorenstein,
then} $\lid_{\Lambda}\Lambda=k$ {\it if and only if} $\lfindim
\Lambda =k$, {\it and} $\rid_{\Lambda}\Lambda=k$ {\it if and only
if} $\rfindim \Lambda =k$.

\vspace{0.2cm}

{\it Proof.} It is well known that $\lid_{\Lambda}\Lambda
\geq\lfindim \Lambda$ and $\rid_{\Lambda}\Lambda \geq\rfindim
\Lambda$.

Assume that $\lfd_{\Lambda}I'_i\leq i+1$ for any $0\leq i \leq k$.
Then, by the dual statement of [4, Theorem 4.7], we have that
$\grade {\rm Ext}_{\Lambda}^i (M, \Lambda)\geq i$ and $\sgrade {\rm
Ext}_{\Lambda}^{i+1} (N, \Lambda)\geq i$ for any $M\in\mod \Lambda$,
$N\in\mod \Lambda ^{op}$ and $1 \leq i \leq k+1$. If $\lfindim
\Lambda =k$, then Ext$_{\Lambda}^{k+1} (M, \Lambda)=0$ by Theorem
2.7, hence $\lid_{\Lambda}\Lambda \leq k$. On the other hand, if
$\rfindim \Lambda =k$, then Ext$_{\Lambda}^{k+2} (M, \Lambda)=0$ by
the dual statement of Theorem 2.7, hence $\rid_{\Lambda}\Lambda \leq
k+1$.

In particular, if $\Lambda$ is $(k+1)$-Gorenstein, then, by
Auslander's Theorem and Theorem 2.7, we get our conclusion
similarly. $\blacksquare$

\vspace{0.2cm}

{\bf Corollary 2.15} {\it If} $\lfd_{\Lambda}I'_i\leq i+1$ {\it for
any} $i \geq 0$, {\it then} $\lid_{\Lambda}\Lambda=\lfindim \Lambda$
{\it and} $\rfindim \Lambda \leq\rid_{\Lambda}\Lambda \leq\rfindim
\Lambda +1$. {\it In particular, if} $\Lambda$ {\it is an Auslander
ring, then} $\lid_{\Lambda}\Lambda=\lfindim \Lambda$ {\it and}
$\rid_{\Lambda}\Lambda=\rfindim \Lambda$.

\vspace{0.5cm}

\centerline{\large \bf 3. Flat dimension and grade of modules}

\vspace{0.2cm}

Recall that $\Lambda$ has {\it dominant dimension} at least $k$,
written $\domdim \Lambda \geq k$, if each $I_i$ is flat for any
$0\leq i \leq k-1$. We write $\domdim \Lambda =\infty$ if $I_i$ is
flat for all $i$. In addition, we denote $K_i=\Ker (I_i \to
I_{i+1})$ for any $i\geq 0$.

\vspace{0.2cm}

{\bf Proposition 3.1} $\domdim \Lambda =\infty$ {\it if and only if}
$\sgrade {\rm Ext}_{\Lambda}^1(M, \Lambda)=\infty$ {\it for any}
$M\in\mod \Lambda$.

\vspace{0.2cm}

{\it Proof.} {\it The sufficiency.} Assume that $\sgrade {\rm
Ext}_{\Lambda}^1(M, \Lambda)=\infty$ for any $M\in\mod \Lambda$.

We will prove that ${\rm Hom}_{\Lambda}({\rm Ext}_{\Lambda}^1(M,
\Lambda), I_i)=0$ for any $M\in\mod \Lambda$ and $i\geq 0$ by using
induction on $i$. We first claim that ${\rm Hom}_{\Lambda}({\rm
Ext}_{\Lambda}^1(M, \Lambda), I_0)=0$. Otherwise, there is a
non-zero homomorphism $f: {\rm Ext}_{\Lambda}^1(M, \Lambda) \to
I_0$. Then $\Im f\bigcap \Lambda \neq 0$ since $\Lambda$ is
essential in $I_0$. So there is a submodule $f^{-1}(\Im f\bigcap
\Lambda)$ of ${\rm Ext}_{\Lambda}^1(M, \Lambda)$ such that
Hom$_{\Lambda} (f^{-1}(\Im f\bigcap \Lambda), \Lambda)\neq 0$, which
contradicts that $\sgrade {\rm Ext}_{\Lambda}^1(M, \Lambda)=\infty$.
Thus we conclude that ${\rm Hom}_{\Lambda}({\rm Ext}_{\Lambda}^1(M,
\Lambda), I_0)=0$. Now suppose $i\geq 1$. Consider the exact
sequence
$$0\to {\rm Hom}_{\Lambda}({\rm
Ext}_{\Lambda}^1(M, \Lambda), K_{i-1}) \to {\rm Hom}_{\Lambda}({\rm
Ext}_{\Lambda}^1(M, \Lambda), I_{i-1})$$ $$\to {\rm
Hom}_{\Lambda}({\rm Ext}_{\Lambda}^1(M, \Lambda), K_i) \to {\rm
Ext}_{\Lambda}^i({\rm Ext}_{\Lambda}^1(M, \Lambda), \Lambda) \to 0$$
for any $i\geq 1$. Since $\sgrade {\rm Ext}_{\Lambda}^1(M,
\Lambda)=\infty$, ${\rm Ext}_{\Lambda}^i({\rm Ext}_{\Lambda}^1(M,
\Lambda), \Lambda)=0$. On the other hand, by induction hypothesis we
have ${\rm Hom}_{\Lambda}({\rm Ext}_{\Lambda}^1(M, \Lambda),
I_{i-1})=0$. So, by the above exact sequence, we have ${\rm
Hom}_{\Lambda}({\rm Ext}_{\Lambda}^1(M, \Lambda), K_i)=0$. By using
an argument similar to the proof of the case $i=0$ we then get that
${\rm Hom}_{\Lambda}({\rm Ext}_{\Lambda}^1(M, \Lambda), I_i)=0$. The
assertion is proved.

By [8, Chapter VI, Proposition 5.3], we have that
Tor$_1^{\Lambda}(I_i, M)\cong {\rm Hom}_{\Lambda}({\rm
Ext}_{\Lambda}^1(M, \Lambda), I_i)=0$ for any $M\in\mod \Lambda$,
which implies that $I_i$ is flat for any $i\geq 0$.

{\it The necessity.} If $\domdim \Lambda =\infty$, then $\Lambda$ is
an Auslander ring and $\sgrade {\rm Ext}_{\Lambda}^1(M, \Lambda)\geq
1$ for any $M\in\mod \Lambda$. Let $X$ be a submodule of
Ext$_{\Lambda}^1(M, \Lambda)$. Then $\grade X\geq 1$. Consider the
exact sequence
$$0\to {\rm Hom}_{\Lambda}(X, K_i) \to
{\rm Hom}_{\Lambda}(X, I_i) \to {\rm Hom}_{\Lambda}(X, K_{i+1}) \to
{\rm Ext}_{\Lambda}^{i+1}(X, \Lambda) \to 0$$ for any $i\geq 0$.
Since $\domdim \Lambda =\infty$, $I_i$ is flat for any $i\geq 0$. By
[6, Theorem 5.2.7], each $I_i$ is a direct limit of free modules in
$\mod \Lambda ^{op}$. So ${\rm Hom}_{\Lambda}(X, I_i)=0$ and hence
${\rm Hom}_{\Lambda}(X, K_i)=0$. By the exactness of the above
sequence we have ${\rm Ext}_{\Lambda}^i(X, \Lambda)=0$ for any
$i\geq 0$ and $\grade X=\infty$. $\blacksquare$

\vspace{0.2cm}

The generalized Nakayama conjecture has an equivalent version as
follows: $\grade S<\infty$ for any simple module $S$ in $\mod
\Lambda ^{op}$ over an artin algebra $\Lambda$ (see [3]); and the
strong Nakayama conjecture says that $\grade N=\infty$ implies $N=0$
for any $N$ in $\mod \Lambda ^{op}$ (see [9]). It is clear that the
generalized Nakayama conjecture is a special case of the strong
Nakayama conjecture. The following result shows that these
conjectures are equivalent if the right flat dimension of each $I_i$
is finite.

\vspace{0.2cm}

{\bf Proposition 3.2} {\it If} $\rfd_{\Lambda}I_i <\infty$ {\it for
all} $i$, {\it then the following statements are equivalent.}

(1) $\grade N=\infty$ {\it implies} $N=0$ {\it for any} $N$ {\it in}
$\mod \Lambda ^{op}$.

(2) $\grade S<\infty$ {\it for any simple module} $S$ {\it in} $\mod
\Lambda ^{op}$.

(3) $\bigoplus _{i\geq 0}I_i$ {\it is an injective cogenerator for
the category of right} $\Lambda$-{\it modules}.

\vspace{0.2cm}

{\it Proof.} $(1)\Rightarrow (2)$ It is trivial.

$(2)\Rightarrow (3)$ Let $S$ be a simple module in $\mod \Lambda
^{op}$. By (2), $\grade S<\infty$ and there is a non-negative
integer $t$ such that Ext$_{\Lambda}^t(S, \Lambda)\neq 0$.

Consider the exact sequences
$$0\to {\rm Hom}_{\Lambda}(S, K_i) \to
{\rm Hom}_{\Lambda}(S, I_i)$$ and $${\rm Hom}_{\Lambda}(S, K_i) \to
{\rm Ext}_{\Lambda}^i(S, \Lambda)\to 0,$$ where $K_i=\Ker (I_i \to
I_{i+1})$ for any $i\geq 0$. Since ${\rm Ext}_{\Lambda}^t(S,
\Lambda)\neq 0$, ${\rm Hom}_{\Lambda}(S, K_t)\neq 0$. So ${\rm
Hom}_{\Lambda}(S, I_t)\neq 0$ and ${\rm Hom}_{\Lambda}(S, \bigoplus
_{i\geq 0}I_i)\neq 0$. It then follows from [1, Proposition 18.15]
that $\bigoplus _{i\geq 0}I_i$ is an injective cogenerator for the
category of right $\Lambda$-modules.

$(3)\Rightarrow (1)$ Let $N \in\mod \Lambda ^{op}$ with $\grade
N=\infty$ and
$$\cdots \to Q_i\to \cdots \to Q_1 \to Q_0 \to N \to 0$$ a
projective resolution of $N$ in $\mod \Lambda ^{op}$. Put
$X_i=\Coker (Q_{i-1}^* \to Q_i^*)$ for any $i\geq 1$. By Lemma 2.8
we have $N\cong$Ext$_{\Lambda}^i(X_i, \Lambda)$ for any $i\geq 1$.

Without loss of generality, we assume that
$\rfd_{\Lambda}I_i=n_i(<\infty)$ for any $i\geq 0$. It follows from
[8, Chapter VI, Proposition 5.3] that ${\rm Hom}_{\Lambda}(N,
I_i)\cong {\rm Hom}_{\Lambda}({\rm Ext}_{\Lambda}^{n_i+1}(X_{n_i+1},
\Lambda), I_i)\cong {\rm Tor}_{n_i+1}^{\Lambda}(I_i, X_{n_i+1})$
\linebreak $=0$ for any $i\geq 0$. So ${\rm Hom}_{\Lambda}(N,
\bigoplus _{i\geq 0}I_i)=0$. However, $\bigoplus _{i\geq 0}I_i$ is
an injective cogenerator for the category of right
$\Lambda$-modules. We then conclude that $N=0$. $\blacksquare$

\vspace{0.2cm}

The famous Nakayama conjecture says that an artin algebra $\Lambda$
is self-injective if $\domdim \Lambda =\infty$. From [3] we know
that the generalized Nakayama conjecture implies the Nakayama
conjecture. By Propositions 3.1 and 3.2 we give here a simple proof
of this implication. Assume that the generalized Nakayama conjecture
is true. If $\domdim \Lambda =\infty$, then $\sgrade {\rm
Ext}_{\Lambda}^1(M, \Lambda)=\infty$ for any $M\in\mod \Lambda$ by
Proposition 3.1. It follows from Proposition 3.2 that
Ext$_{\Lambda}^1(M, \Lambda)=0$ for any $M\in\mod \Lambda$ and
$\Lambda$ is self-injective.

In the rest of this section, we study the properties of pure
modules. Let $k$ be a non-negative integer. A non-zero module $M$ in
$\mod \Lambda$ (resp. $\mod \Lambda ^{op}$) is said to be {\it pure}
of grade $k$ if $\grade A=k$ for each non-zero submodule $A$ of $M$.
The notion of pure modules here coincides with that given in [7]
when $\Lambda$ is Auslander-Gorenstein.

\vspace{0.2cm}

{\bf Lemma 3.3}  {\it Let} $\Lambda$ {\it be a 1-Gorenstein ring}.
{\it Then a module} $N$ {\it in} $\mod \Lambda ^{op}$ {\it is pure
of grade 0 if and only if it is torsionless}.

\vspace{0.2cm}

{\it Proof.} If $N$ is torsionless, then each non-zero submodule $A$
of $N$ is also torsionless and so $A^*\neq 0$, that is, $\grade
A=0$. Conversely, assume that $N$ is pure of grade 0. By Lemma 2.8
there is a module $A\in\mod \Lambda$ such that $\Ker \sigma _N
\cong$Ext$_{\Lambda}^1(M, \Lambda)$. Since $\Lambda$ is
1-Gorenstein,  $[{\rm Ext}_{\Lambda}^1(M, \Lambda)]^*=0$. So
Ext$_{\Lambda}^1(M, \Lambda)=0$ since $N$ is pure of grade 0. It
follows that $\sigma _N$ is a monomorphism and $N$ is torsionless.
$\blacksquare$

\vspace{0.2cm}

{\it Remark.} The proof of Lemma 3.3 in fact proves the following
more general result. If $[{\rm Ext}_{\Lambda}^1(M, \Lambda)]^*$
\linebreak $=0$ for any $M\in\mod \Lambda$, then a module $N$ in
$\mod \Lambda ^{op}$ is pure of grade 0 if and only if it is
torsionless. So, by the dual statement of [4, Theorem 4.7] we have
that if $\Lambda$ is a noetherian algebra with $\lfd_{\Lambda}I'_0
\leq 1$, then a module $N$ in $\mod \Lambda ^{op}$ is pure of grade
0 if and only if it is torsionless.

\vspace{0.2cm}

{\bf Lemma 3.4} {\it Let} $k$ {\it be a positive integer and} $m
\geq -1$ {\it an integer. Then the following statements are
equivalent.}

(1) $\rfd_{\Lambda}\bigoplus _{i=0}^{k-1}I_i \leq k+m$.

(2) $\sgrade {\rm Ext}_{\Lambda}^{k+m+1}(M, \Lambda)\geq k$ {\it for
any} $M \in\mod \Lambda$.

\vspace{0.2cm}

{\it Proof.} It was proved in [13, Theorem 2.8] in case $m$ is a
non-negative integer. When $m=-1$, the proof is similar to that of
[13, Theorem 2.8], we omit it. $\blacksquare$

\vspace{0.2cm}

The following lemma is a dual statement of Lemma 3.4.

\vspace{0.2cm}

{\bf Lemma 3.5} {\it Let} $k$ {\it be a positive integer and} $m
\geq -1$ {\it an integer. Then the following statements are
equivalent.}

(1) $\lfd_{\Lambda}\bigoplus _{i=0}^{k-1}I'_i \leq k+m$.

(2) $\sgrade {\rm Ext}_{\Lambda}^{k+m+1}(N, \Lambda)\geq k$ {\it for
any} $N \in\mod \Lambda ^{op}$.

\vspace{0.2cm}

We are now in a position to give the main result in this section.

\vspace{0.2cm}

{\bf Theorem 3.6} {\it Let} $k$ {\it be a non-negative integer and}
$N$ {\it in} $\mod \Lambda ^{op}$ {\it with} $\grade N=k<\infty$.

(1) {\it If} $\lfd_{\Lambda}\bigoplus _{i=0}^k I'_i \leq k$ {\it
and} $N$ {\it is pure of grade} $k$, {\it then} $N$ {\it can be
embedded into a finite direct sum of copies of} $I_k$.

(2) {\it If} $\rfd_{\Lambda}\bigoplus _{i=0}^{k-1}I_i \leq k-1$,
$\rfd_{\Lambda}I_k \leq k$ {\it and} $N$ {\it can be embedded into a
finite direct sum of copies of} $I_k$, {\it then} $N$ {\it is pure
of grade} $k$.

\vspace{0.2cm}

{\it Proof.} The case $k=0$ follows from Lemma 3.3. Now suppose
$k\geq 1$.

$(1)$ Also put $K_i=\Ker (I_i \to I_{i+1})$ for any $i\geq 0$. Since
$N$ is pure of grade $k$, $\grade X=k$ for any submodule $X$ of $N$.
Then it is not difficult to see that Hom$_{\Lambda}(N, I_i)=0$ for
any $0\leq i \leq k-1$ and so Ext$_{\Lambda}^k(N,
\Lambda)\cong$Hom$_{\Lambda}(N, K_k)$.

Let $\eta _1$, $\eta _2$, $\cdots$, $\eta _n$ be a set of generators
of Hom$_{\Lambda}(N, K_k)$ in End$_{\Lambda}(K_k)$. Put $\eta =(\eta
_1, \eta _2, \cdots, \eta _n)' : N\to K_k^{(n)}$, $U=\Ker \eta$ and
$V=\Im \eta$. We use $\pi : N\to V$ to denote the natural
epimorphism. Since $\grade N=k$ and $\grade U\geq k$ by assumption,
$\grade V\geq k$. By using an argument similar to the above we then
have Ext$_{\Lambda}^k(V, \Lambda)\cong$Hom$_{\Lambda}(V, K_k)$. In
addition, it is not difficult to see that Hom$_{\Lambda}(\eta, K_k)$
is an epimorphism, so Hom$_{\Lambda}(\pi, K_k): {\rm
Hom}_{\Lambda}(V, K_k)\to {\rm Hom}_{\Lambda}(N, K_k)$ is also an
epimorphism and hence an isomorphism.

On the other hand, we have an exact sequence:
$$0={\rm Ext}_{\Lambda}^{k-1}(U, \Lambda)\to
{\rm Ext}_{\Lambda}^k(V, \Lambda)\buildrel {{\rm
Ext}_{\Lambda}^k(\pi, \Lambda)} \over \longrightarrow {\rm
Ext}_{\Lambda}^k(N, \Lambda)\to {\rm Ext}_{\Lambda}^k(U, \Lambda)\to
{\rm Ext}_{\Lambda}^{k+1}(V, \Lambda).$$ So ${\rm
Ext}_{\Lambda}^k(\pi, \Lambda)$ is an isomorphism and ${\rm
Ext}_{\Lambda}^k(U, \Lambda)$ is isomorphic to a submodule of ${\rm
Ext}_{\Lambda}^{k+1}(V, \Lambda)$. Since $\lfd_{\Lambda}\bigoplus
_{i=0}^k I'_i \leq k$, $\grade {\rm Ext}_{\Lambda}^k(U, \Lambda)\geq
k+1$ by Lemma 3.5. It then follows from Lemma 2.3 that ${\rm
Ext}_{\Lambda}^k(U, \Lambda)=0$ and $\grade U\geq k+1$, which
implies that $U=0$ and $\eta$ is a monomorphism.

$(2)$  Since $\grade N=k$, by Lemma 2.8 there is an $X\in\mod
\Lambda$ such that $N\cong$Ext$_{\Lambda}^k(X, \Lambda)$. On the
other hand, notice that $\rfd_{\Lambda}\bigoplus _{i=0}^{k-1}I_i
\leq k-1$, so $\sgrade N=\sgrade {\rm Ext}_{\Lambda}^k(X,
\Lambda)\geq k$ by Lemma 3.4.

Let $U$ be a submodule of $N$. Then $\grade U\geq k$. If $\grade
U\geq k+1$, then, again by Lemma 2.8, there is a $Y\in\mod \Lambda$
such that $U\cong$Ext$_{\Lambda}^{k+1}(Y, \Lambda)$. It follows from
[8, Chapter VI, Proposition 5.3] that ${\rm Hom}_{\Lambda}(U,
I_k)\cong {\rm Hom}_{\Lambda}({\rm Ext}_{\Lambda}^{k+1}(Y, \Lambda),
I_k)\cong {\rm Tor}_{k+1}^{\Lambda}(I_k, Y)=0$ since
$\rfd_{\Lambda}I_k \leq k$. However, $U$ can be embedded into a
finite direct sum of copies of $I_k$, so $U=0$. This completes the
proof. $\blacksquare$

\vspace{0.2cm}

Especially, we have the following.

\vspace{0.2cm}

{\bf Theorem 3.7} {\it Let} $k$ {\it be a non-negative integer and}
$\Lambda$ {\it a} $(k+1)$-{\it Gorenstein ring. Then the following
statements are equivalent for a module} $N\in\mod \Lambda ^{op}$
{\it with} $\grade N=k<\infty$.

(1) $N$ {\it is pure of grade} $k$.

(2) $N$ {\it can be embedded into a finite direct sum of copies
of} $I_k$.

\vspace{0.2cm}

{\bf Corollary 3.8} ([12, Theorem 6.3]) {\it Let} $I$ {\it be an
indecomposable injective right} $\Lambda$-{\it module with}
$\rfd_{\Lambda}I =k<\infty$. {\it If} $\rfd_{\Lambda}\bigoplus
_{i=0}^{k-1}I_i \leq k-1$ {\it and} $\lfd_{\Lambda}\bigoplus
_{i=0}^k I'_i \leq k$, {\it then} $I$ {\it appears as a direct
summand of} $I_k$.

\vspace{0.2cm}

{\it Proof.} We first prove the case $k=0$. Let $0\neq X$ be a
finitely generated submodule of $I$. Then $I$ is the injective
envelope of $X$. By [18, Theorem 1.2], $X$ is torsionless. So $X$
can be embedded into a finite direct sum of copies of $\Lambda$ and
hence $I$ can be embedded into a finite direct sum of copies of
$I_0$, which yields that $I$ is isomorphic to a direct summand of
$I_0$.

Now suppose $k\geq 1$. Put $E=\bigoplus _{i=0}^{k-1}I_i$. Then
$\rfd_{\Lambda}E\leq k-1$. Notice that $\rfd_{\Lambda}I= k$, so $E$
does not cogenerate $I$ and hence there is a submodule $X$ of $I$
such that Hom$_{\Lambda}(X, E)=0$. We may assume that $X$ is
finitely generated. Let $Y$ be a submodule of $X$. Clearly, we have
Hom$_{\Lambda}(Y, E)=0$. It follows that $\grade Y\geq k$. If
$\grade Y\geq k+1$, then by Lemma 2.8, there is a module $M\in\mod
\Lambda$ such that $Y\cong$Ext$_{\Lambda}^{k+1}(M, \Lambda)$. Thus,
by [8, Chapter VI, Proposition 5.3], we have ${\rm Hom}_{\Lambda}(Y,
I)\cong {\rm Hom}_{\Lambda}({\rm Ext}_{\Lambda}^{k+1}(M, \Lambda),
I)\cong {\rm Tor}_{k+1}^{\Lambda}(I, M)=0$ since
$\rfd_{\Lambda}I=k$, which is a contradiction. Therefore $X$ is pure
of grade $k$. By Theorem 3.6, $X$, and hence $I$, can be embedded
into a finite direct sum of copies of $I_k$. This completes the
proof. $\blacksquare$

\vspace{0.2cm}

Recall that $\Lambda$ is called an {\it Auslander-Gorenstein ring}
if it is an Auslander ring with finite left and right self-injective
dimensions.

\vspace{0.2cm}

{\bf Corollary 3.9} {\it Let} $\Lambda$ {\it be an
Auslander-Gorenstein ring and} $M$ {\it in} $\mod \Lambda$ {\it
with} $\grade M=k$. {\it Then} Ext$_{\Lambda}^k(M, \Lambda)$ {\it
can be embedded into a finite direct sum of} $I_k$.

\vspace{0.2cm}

{\it Proof.} By [15, Corollary 3.7], Ext$_{\Lambda}^k(M, \Lambda)$
is pure of grade $k$. It follows from Theorem 3.7 that
Ext$_{\Lambda}^k(M, \Lambda)$ can be embedded into a finite direct
sum of $I_k$. $\blacksquare$

\vspace{0.2cm}

{\bf Proposition 3.10} {\it Let} $\Lambda$ {\it be an
Auslander-Gorenstein ring with}
$\lid_{\Lambda}\Lambda=\rid_{\Lambda}\Lambda=k$. {\it Then}
Ext$_{\Lambda}^k(M, \Lambda)$ {\it is pure of grade} $k$ {\it for
any} $M\in\mod \Lambda$ {\it with} Ext$_{\Lambda}^k(M, \Lambda)\neq
0$.

\vspace{0.2cm}

{\it Proof.} Let $M$ be in $\mod \Lambda$ with Ext$_{\Lambda}^k(M,
\Lambda)\neq 0$. Since $\Lambda$ is an Auslander ring, $\sgrade {\rm
Ext}_{\Lambda}^k(M, \Lambda) \geq k$. Let $Y$ be a submodule of
Ext$_{\Lambda}^k(M, \Lambda)$ in $\mod \Lambda ^{op}$. Then $\grade
Y\geq k$. If $\grade Y\geq k+1$, then $\grade Y=\infty$ since
$\rid_{\Lambda}\Lambda=k$. It follows from [9, Theorem 2] that
$Y=0$. This completes the proof. $\blacksquare$

\vspace{0.2cm}

By [16, Proposition 1], we have that if
$\lid_{\Lambda}\Lambda=\rid_{\Lambda}\Lambda=k$ then
$\rfd_{\Lambda}I_k=k$. So a $k$-Gorenstein ring with
$\lid_{\Lambda}\Lambda=\rid_{\Lambda}\Lambda=k$ is just an
Auslander-Gorenstein ring with
$\lid_{\Lambda}\Lambda=\rid_{\Lambda}\Lambda=k$. By Theorem 3.7 and
Proposition 3.10 we then get the following.

\vspace{0.2cm}

{\bf Corollary 3.11} ([16, Corollary 5]) {\it Let} $\Lambda$ {\it be
a} $k$-{\it Gorenstein ring with}
$\lid_{\Lambda}\Lambda=\rid_{\Lambda}\Lambda=k$ {\it and} $M$ {\it
in} $\mod \Lambda$. {\it Then} Ext$_{\Lambda}^k(M, \Lambda)$ {\it
can be embedded into a finite direct sum of} $I_k$. {\it Moreover,
if} Ext$_{\Lambda}^k(M, \Lambda)\neq 0$, {\it then}
Ext$_{\Lambda}^k(M, \Lambda)$ {\it is pure of grade} $k$.

\vspace{0.2cm}

Recall from [1] that a non-zero right $\Lambda$-module $H$ is called
{\it uniform} if each of its non-zero submodules is essential in
$H$.

\vspace{0.2cm}

{\bf Corollary 3.12} {\it Let} $\Lambda$ {\it be an
Auslander-Gorenstein ring with}
$\lid_{\Lambda}\Lambda=\rid_{\Lambda}\Lambda=k$ {\it and} $N$ {\it a
uniform module in} $\mod \Lambda ^{op}$. {\it If} $N$ {\it is pure
of grade} $k$, {\it then} $N$ {\it can be embedded into} $I_k$, {\it
but cannot be embedded into} $\bigoplus _{i=0}^{k-1}I_i$ {\it (if}
$k\geq 1${\it )}.

\vspace{0.2cm}

{\it Proof.} Let $N$ be a uniform module in $\mod \Lambda ^{op}$.
Then $E(N)$ (the envelope of $N$) is indecomposable. If $N$ is pure
of grade $k$, then, by Theorem 3.7, $N$ and $E(N)$ can be embedded
into a finite direct sum of $I_k$. So $E(N)$ is isomorphic to a
direct summand of $I_k$ and hence $N$ can be embedded into $I_k$. On
the other hand, if $k\geq 1$, then, by [16, Corollary 7], we have
that $I_k$ and $\bigoplus _{i=0}^{k-1}I_i$ have no isomorphic direct
summands in common. It follows from the above argument that $E(N)$,
and hence $N$, cannot be embedded into $\bigoplus _{i=0}^{k-1}I_i$.
$\blacksquare$

\vspace{0.5cm}

\centerline{\large \bf 4. The socle of the last term in a
minimal injective resolution}

\vspace{0.2cm}

For a right $\Lambda$-module $X$, the unique largest semisimple
submodule of $X$ is called the {\it socle} of $X$, and denoted by
Soc$(X)$ (see [1]). In this section we show that under some grade
conditions of modules the socle of $I_k$ is non-zero. In fact we
will prove the following.

\vspace{0.2cm}

{\bf Theorem 4.1} {\it Assume that}
$\lid_{\Lambda}\Lambda=\rid_{\Lambda}\Lambda=k$. {\it If} $\grade
{\rm Ext}_{\Lambda}^k(M, \Lambda)\geq k$ {\it for any} $M\in\mod
\Lambda$ {\it and} $\grade {\rm Ext}_{\Lambda}^i(N, \Lambda)\geq i$
{\it for any} $N\in\mod \Lambda ^{op}$ {\it with} $1\leq i \leq
k-1$, {\it then} Soc$(I_k)\neq 0$.

\vspace{0.2cm}

{\it Proof.} The case $k \leq 2$ was proved in [12, Theorem 4.5].
Now suppose $k\geq 3$.

Since $\lid_{\Lambda}\Lambda=k$, there is a module $M\in\mod
\Lambda$ such that Ext$_{\Lambda}^k(M, \Lambda)\neq 0$. Let $N$ be a
maximal submodule of Ext$_{\Lambda}^k(M, \Lambda)$ and
$$0 \to N \to {\rm Ext}_{\Lambda}^k(M, \Lambda) \to S \to 0$$
an exact sequence in $\mod \Lambda ^{op}$. Then $S$ is simple.

Since $\grade {\rm Ext}_{\Lambda}^i(N, \Lambda)\geq i$ for any
$N\in\mod \Lambda ^{op}$ and $1\leq i \leq k-1$ by assumption,
$\rfd_{\Lambda}I_i \leq i+1$ for any $0 \leq i \leq k-2$ by [4,
Theorem 0.1] and the remark following it. It follows from [8,
Chapter VI, Proposition 5.3] that Hom$_{\Lambda}({\rm
Ext}_{\Lambda}^k(M, \Lambda), I_i)\cong$Tor$_k^{\Lambda}(I_i, M)=0$
for any $0 \leq i \leq k-2$. We then have Hom$_{\Lambda}(N,
I_i)=0=$Hom$_{\Lambda}(S, I_i)$ for any $0 \leq i \leq k-2$, and
therefore $\grade N\geq k-1$ and $\grade S\geq k-1$.

On the other hand, $\grade {\rm Ext}_{\Lambda}^k(M, \Lambda)\geq k$
by assumption, so we have Ext$_{\Lambda}^{k-1}(S,
\Lambda)\cong$Ext$_{\Lambda}^{k-2}(N, \Lambda)=0$ and hence $\grade
S\geq k$. If $\grade S\geq k+1$, then $\grade S=\infty$ since
$\rid_{\Lambda}\Lambda=k$. It follows from Lemma 2.8 that $S$ is
reflexive and $S=0$, which is a contradiction. So we conclude that
$\grade S=k$ and Ext$_{\Lambda}^k(S, \Lambda)\neq 0$. Thus we have
Hom$_{\Lambda}(S, I_k)\neq 0$, which implies that $S$ is isomorphic
to a simple submodule of $I_k$ and Soc$(I_k)\neq 0$. $\blacksquare$

\vspace{0.2cm}

{\bf Corollary 4.2} {\it Assume that} $\Lambda$ {\it is a noetherian
algebra with} $\lid_{\Lambda}\Lambda=\rid_{\Lambda}\Lambda=k$. {\it
If} $\rfd_{\Lambda}I_i \leq i+1$ {\it for any} $0\leq i \leq k-2$
{\it and} $\lfd_{\Lambda}I'_i \leq i+1$ {\it for any} $0\leq i \leq
k-1$, {\it then} Soc$(I_k)\neq 0$.

\vspace{0.2cm}

{\it Proof.} Since $\rfd_{\Lambda}I_i \leq i+1$ for any $0\leq i
\leq k-2$, by [4, Theorem 4.7] we have $\grade {\rm
Ext}_{\Lambda}^i(N, \Lambda)\geq i$ for any $N\in\mod \Lambda ^{op}$
and $1\leq i \leq k-1$. On the other hand, since $\lfd_{\Lambda}I'_i
\leq i+1$ for any $0\leq i \leq k-1$, by the dual statement of [4,
Theorem 4.7] we have $\grade {\rm Ext}_{\Lambda}^i(M, \Lambda)\geq
i$ for any $M\in\mod \Lambda$ and $1\leq i \leq k$. We now get our
conclusion by Theorem 4.1. $\blacksquare$

\vspace{0.2cm}

From Auslander's Theorem and Theorem 4.1 we obtain the following.

\vspace{0.2cm}

{\bf Corollary 4.3} ([11, Proposition 1.1]) {\it Let} $\Lambda$ {\it
be an Auslander-Gorenstein ring with}
$\lid_{\Lambda}\Lambda=\rid_{\Lambda}\Lambda=k$. {\it Then}
Soc$(I_k)\neq 0$.

\vspace{0.2cm}

{\bf Example 4.4} There are rings satisfying the assumption in
Theorem 4.1 but not satisfying the assumption in Corollary 4.3. Let
$K$ be a field and $\Lambda$ the finite dimensional path algebra
over $K$ given by the quiver

$$
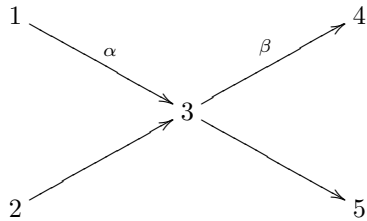
\xymatrix{1\ar[rrd]^\alpha &&&&4\\
&&3\ar[rru]^\beta\ar[rrd]\\
2\ar[rru]&&&&5}$$ modulo the ideal generated by $\beta \alpha$. Then
$\lfd_{\Lambda}I'_0$=
$\lfd_{\Lambda}I'_1=\rfd_{\Lambda}I_0=\rfd_{\Lambda}I_1$=1,
$\lfd_{\Lambda}I'_2=\rfd_{\Lambda}I_2$=2 and
$\lid_{\Lambda}\Lambda=\rid_{\Lambda}\Lambda$=2. By [4, Theorem 4.7]
and its dual statement, for any $i\geq 1$ we have $\grade {\rm
Ext}_{\Lambda}^i(M, \Lambda)\geq i$ for any $M\in\mod \Lambda$ and
$\grade {\rm Ext}_{\Lambda}^i(N, \Lambda)\geq i$ for any $N\in\mod
\Lambda ^{op}$. So $\Lambda$ satisfies the assumption in Theorem 4.1
(In fact $\Lambda$ also satisfies the assumption in Corollary 4.2).
But $\Lambda$ is clearly not Auslander-Gorenstein.

\vspace{0.2cm}

{\bf Acknowledgements} The research of the author was partially
supported by the Specialized Research Fund for the Doctoral Program
of Higher Education (Grant No. 20060284002) and NSF of Jiangsu
Province of China (Grant No. BK2007517).

\end{document}